\documentclass{amsart}

\usepackage[left=3cm,top=3cm,right=3cm,bottom=3cm]{geometry}
\usepackage{amssymb, amsmath, amsthm, units}
\usepackage{verbatim}
\usepackage{comment}
\usepackage{graphicx}
\usepackage{enumerate}
\usepackage{color}

\DeclareGraphicsExtensions{.pdf,.jpg,.png}

\theoremstyle{plain}
\newtheorem{theorem}{Theorem}
\newtheorem{lemma}[theorem]{Lemma}

\newtheorem{corollary}[theorem]{Corollary}

\newtheorem{problem}[theorem]{Problem}

\theoremstyle{definition}

\newcommand{\Z}{\mathbb{Z}}
\newcommand{\N}{\mathbb{N}}
\newcommand{\F}{\mathbb{F}}
\newcommand{\M}{\mathcal{M}}
\newcommand{\R}{\mathbb{R}}

\makeatletter

\begin{document}

\title{The sum of nonsingular matrices is often nonsingular}

\author{J\'{o}zsef Solymosi}
\address{\noindent Department of Mathematics, 
University of British Columbia, 1984 Mathematics Road, 
Vancouver, BC, V6T 1Z2, Canada}
\email{solymosi@math.ubc.ca}
\thanks{Research supported in part by a  NSERC and an OTKA NK 104183 grant}

\date{}

\begin{abstract} If $\M$ is a set of nonsingular $k\times k$ matrices then for many pairs of matrices, $A,B\in\M,$ the sum is nonsingular, $\det(A+B)\neq 0.$ We prove a more general statement on nonsingular sums with a geometric application. 
\end{abstract}

\maketitle
\noindent

\medskip
\section{Introduction}
It is a simple fact in linear algebra, that while the product of nonsingular matrices is nonsingular, the similar statement is false for the sum. On the other hand, one can expect some control over such sums since most matrices are nonsingular. (The proper notion of {\em most matrices} depends on the underlying field. Much more details on this subject can be find in a book of Terry Tao \cite{TAO}.) Finding the inverse or a generalized inverse of the sum of two matrices has important applications in mathematics and in applications. For the review and the history on deriving the inverse of the sum of matrices we refer to \cite{HENDERSON} and \cite{BT}. In this paper we show that---under some mild conditions---if $\M$ is a set of  $k\times k$ matrices then for many pairs, $A,B\in\M,$ the sum is nonsingular. Our main tool is the so-called ``Polynomial Method"\footnote{ In many cases the {\em Linear Algebra Method} or {\em Rank Method} would be a better description.} which has been used in combinatorics since the 70's and has proven to be very useful in a number of problems. There is a nice lecture book by Larry Guth reviewing old and new applications of the method \cite{GUTH}. There is a striking, very recent application to the cap-set problem by Croot, Lev, and Pach \cite{CLP} and Ellenberg and Gijswijt \cite{EG}. The latter two inspired a large number of interesting results using a counting method similar to what we will apply in this work. 

Over finite fields Anderson and Badawi investigated the graph where the vertices are the elements of $SL_n(\F_q)$ and two matrices, $A,B\in SL_n(\F_q)$ form an edge if $\det(A+B)\neq 0.$ Akbari, Jamaali and Fakhari \cite{AJF} proved that the clique number of such graphs is bounded by a universal constant, independent of $\F_q$ for odd $q.$ We will refine their result giving almost sharp bound  on the clique number. Tomon \cite{TO} showed that the chromatic number of this graph is at least $(q/4)^{\lfloor{n/2}\rfloor}.$

\section{Results}
First we state an important case of our main result below. For any set of nonsingular matrices a positive fraction of the pairs add up to a nonsingular matrix. 

\begin{theorem}\label{first}
For every $k\in \N$ there is a constant, $c>0,$ depending on $k$ only, such that if $\M$ is a set of nonsingular $k\times k$ matrices over a field $K,$ characteristic $\neq 2$ then the number of pairs, $A,B\in\M,$ such that $\det(A+B)\neq 0,$ is at least $c|\M|^2.$
\end{theorem}

We postpone the proof until after our next theorem, where we are going to give a necessary and sufficient condition under which, for a positive fraction of the pairs,  $A,B\in\M,$ $\det(A+B)\neq 0.$ For the exact statement we are going to consider two (not necessary disjoint) sets of $k\times k$ matrices, $\M_1,\M_2$ as the two vertex sets of a bipartite graph, $G(\M_1,\M_2),$ where $A\in\M_1$ and $B\in\M_2$ are connected by an edge iff  $\det(A+B)\neq 0.$ In what follows we will suppose that $|\M_1|=|\M_2|=n.$ A {\em matching} in a graph is a set of vertex disjoint edges. It follows from elementary graph theory that if the number of edges in $G(\M_1,\M_2),$ is at least $cn^2$ then it contains a matching of size at least $cn.$ We show that---at least asymptotically---the converse holds as well.

\begin{theorem}\label{main}
For every $k\in \N$ there is a constant, $c>0,$ depending on  $k$ only, such that the following holds: Let $\M_1,\M_2$ be two $n$-element sets of $k\times k$ matrices over a field $K,$ characteristic $\neq 2.$
If the bipartite graph, $G(\M_1,\M_2),$ as defined above, contains a perfect matching ($n$ vertex disjoint edges) then it has at least $cn^2$ edges.
\end{theorem}

Requiring a perfect matching is not a real restriction here. If $G(\M_1,\M_2)$ contains a matching of size $m,$ then one can restrict the graph to the vertices of the matching and applying Theorem \ref{main} guarantees at least $cm^2$ edges. 

\medskip
{\em Proof of Theorem \ref{first}:} Set $\M_1=\M$ and $\M_2=\M.$ Since $\M$ consists of nonsingular matrices, for every $A\in\M$ the pair $(A,A)$ is an edge in $G(\M_1,\M_2).$ These edges form a perfect matching, so applying Theorem \ref{main} we see that the number of pairs, $A,B\in\M,$ such that $\det(A+B)\neq 0,$ is at least $c|\M|^2.$ 
\qed

\medskip
{\em Proof of Theorem \ref{main}:} Let us suppose that $\M_1=\{A_1,A_2,\ldots ,A_n\},$ $\M_2=\{B_1,B_2,\ldots ,B_n\},$ and that the perfect matching in $G(\M_1,\M_2)$ is given by the edges $(A_i,B_i).$ We define an $n\times n$ matrix, $H=(h_{ij}),$ with entries $h_{ij}:=\det(A_i+B_j).$ First we show that $H$ has rank smaller than $4^k.$ We will expand $\det(A_i+B_j)$ into subsums. For that we introduce some notations. For a $k\times k$ matrix, $M=(m_{ij}),$ and two subsets of the index set $I,J\subset[k],$ we define the submatrix with rows from $I$ and columns from $J$, as $M[I\times J]=(m_{ij})_{i\in I, j\in J}.$ The sign of the sums is determined by the $\sigma(I,J)=\sum_{i\in I}i+\sum_{j\in J}j$ function. As usual, $\bar{S}$ denotes the complement of $S.$ in our case $\bar{I}=[k]\setminus I.$

\begin{equation}\label{rank}
\det(A_i+B_j)=\sum_{\ell=0}^k\sum_{\substack{I,J\subset[k],\\ |I|=|J|=\ell}}(-1)^{\sigma(I,J)}\det(A_i[I\times J])\det(B_j[\bar{I}\times\bar{J}]).
\end{equation}

This formula is easy to prove, and it probably has several possible references. There is a nice discussion on the formula with a proof in \cite{MAMA}.
For given subsets $I,J\subset[k],$ $|I|=|J|,$ we define the $n\times n$ matrix $H^{(I,J)}$ as a matrix with entries 

$$h^{(I,J)}_{ij}=(-1)^{\sigma(I,J)}\det(A_i[I\times J])\det(B_j[\bar{I}\times\bar{J}]).$$

Note that $H^{(I,J)}$ has rank at most one, since every column is a multiple of any other column. As $H=\sum H^{(I,J)}$, we have the following upper bound on the rank

\[rank(H)\leq \sum_{i=0}^k{k\choose i}^2={2k\choose k}\sim \frac{4^k}{\sqrt{2k}}.\]

In the second part of the proof we define an auxiliary graph on $n$ vertices, denoted by $G^*_n.$ Two vertices, $v_i$ and $v_j$ are connected by an edge iff $h_{ij}=h_{ji}=0,$ or equivalently, 
\begin{equation}\label{singsum}
\det(A_i+B_j)=\det(A_j+B_i)=0.
\end{equation}
If the number of nonsingular sums between $\M_1$ and $\M_2$ is less than $cn^2,$ then  $G^*_n$ has at least ${n\choose 2}-cn^2\sim \frac{(1-c)n^2}{2}$ edges. On the other hand, as we shall see, this graph can not have more than $\frac{(1-1/4^k)n^2}{2}$ edges. This sets the $4^{-k}$ lower bound on $c.$ Indeed, if the graph, $G^*_n,$  had more 
edges, then by Tur\'an's Theorem \cite{TUR} it would contain a complete subgraph of size $4^k.$ But it would mean that in $H,$ where all entries in the diagonal are non-zero, we find a diagonal submatrix, a leading principal submatrix of size $4^k,$ which contradicts our upper bound on the rank of $H.$

\qed

The rank bounds in the proof of Theorem \ref{main} are exponential in $k.$ It is unavoidable as the following examples show. 

\medskip
\noindent
{\bf Example 1:}
This construction is about a set of  nonsingular $k\times k$ matrices over the reals. Let $\M$ be the collection of the $2^k$ diagonal matrices with $\pm 1$ entries in the diagonal.
The sum of any two distinct matrices is singular. The rank of $H$ in the proof Theorem \ref{main} for this set (with $\M_1=\M_2=\M$) is $2^k.$ This example can be extended to an arbitrary large $\M.$ Instead of $\pm 1$ in the last diagonal entry, it now can be selected from any element of the set $\{-s,-s+1,\ldots, -1,1,2,\ldots, s-1,s\}.$ Then $|\M|=s2^{k},$ all elements are nonsingular and the number of pairs with nonsingular sums is $|\M|(2s-1).$ This example shows that one can not expect better than exponential bound on $c$ in Theorem \ref{first}. However there is still a gap between the proved upper and lower bounds. We know that $4^{-k}<c\leq 2^{-k-1}.$ Maybe the upper bound (the construction) is closer to the truth.

\medskip
\noindent
{\bf Example 2:}
Here we define $\M_1=\{A_1,A_2,\ldots ,A_n\},$ $\M_2=\{B_1,B_2,\ldots ,B_n\},$ with $n$ close to $4^k$ such that $\det(A_i+B_j)=0$ if and only if $i\neq j.$ For the sake of simplicity, let us suppose that $k$ is even. For every possible $I,J\subset[k], |I|=|J|=k/2,$ we embed the identity matrix into the all zero $k\times k$ matrix, to positions $I\times J,$ without changing the order of rows or columns. More formally, if $I=\{i_1,\ldots,i_{k/2}\}$ and $J=\{j_1,\ldots,j_{k/2}\},$ then in $M_{I,J}$ there are zeros  in all positions but the $i_\ell ,j_\ell$ entries, where it is 1 ($1\leq \ell\leq k/2.$) For every such $I,J$ index sets there is one such matrix, $M_{I,J}$. We pair them, if $M_{I,J}=A_r$ is in $\M_1$ then we place $M_{\bar{I},\bar{J}}=B_r$ to $\M_2.$ There are ${k\choose k/2}^2/2\sim 4^k/k$ pairs. Also, $\det(M_{I',J'}+M_{\bar{I},\bar{J}})\neq 0$ iff $I=I'$ and $J=J'$ (In every other case there is an all-zero row or column in the sum). Similar to Example 1, this construction can be extended to arbitrary large sets. Instead of the identity matrix, let's embed diagonal matrices with diagonal entries from the set $S=\{-s,-s+1,\ldots, -1,1,2,\ldots, s-1,s\},$ i.e. embed the $k/2\times k/2$ matrices $tI,$ $t\in S.$ The matrices are $M_{I,J}^t$ and the size of $\M_1$ (same as $|\M_2|$) is $s{k\choose k/2}^2.$ Every matrix 
$M_{I,J}^t\in \M_1$ has exactly $2s$ matrices (of the form $M_{\bar{I},\bar{J}}^{t'}\in\M_2$) such that their sum is nonsingular. This construction shows that the $c\geq 4^{-k}$ bound is almost sharp in Theorem \ref{main}.

\medskip

While the above examples show the limits of possible improvements, there are interesting open questions remain. One of the most important questions is the following.

\begin{problem}
Is it true that for every $k\in \N$ there is a constant, $c>0,$ depending on $k$ only, such that if $\M$ is a set of nonsingular $k\times k$ matrices over $\R,$ then one can always find a subset of the matrices, $\M'\subset\M,$ such that sums in $\M'$ are nonsingular, $\det(A+B)\neq 0$ for any $A,B\in\M',$ and $|\M'|\geq c|\M|$?
\end{problem}

Using standard arguments from Ramsey Theory, it is easy to see that there is always an $\M'\subset\M$ such that the sums in the set are nonsingular and $|\M'|\geq |\M|^c,$ but we expect that a much better bound holds here.  

\section{A geometric application}
In a geometric application we are going to consider $d$-dimensional flats in $\R^{2d}.$ We will prove a type of removal lemma for flats. In preparation, first
we show that if some pairs have a single point intersection, then many intersect in a point only. 

\begin{lemma}\label{flat_pairs}
Let us suppose that we are given $n$ pairs of $d$-dimensional flats in $\R^{2d}$ labelled $F_i,E_i,$ such that $F_i$ and $E_i$ intersect in one point ($1\leq i\leq n$). Then there are at least $n^2/4^d$ pairs of flats, $F_i,E_j,$ which intersect in a single point.
\end{lemma}

{\em Proof:} Let us write the equations of flats as $F_i:\vec{y}=A_i\vec{x} +\vec{v_i}$ and $E_i:\vec{y}=B_i\vec{x} +\vec{w_i},$ for all $1\leq i\leq n.$ Since $F_i$ and $E_i$ intersect in one point, the system of equations $(A_i-B_i)\vec{x}=\vec{v_i}-\vec{w_i}$ has a unique solution, so $A_i-B_i$ is nonsingular for all $1\leq i\leq n.$ We can now apply Theorem \ref{main} to conclude that there are at least $n^2/4^d$ $i,j$ pairs such that $A_i-B_j$ is nonsingular. Then the equation $(A_i-B_j)\vec{x}=\vec{v_i}-\vec{w_j}$ has a unique solution, so $F_i$ and $E_j$ have a single intersection point.
\qed

\medskip

Now we are ready to state and prove our geometric result. Removal lemmas are important tools in graph theory and additive combinatorics. The simplest version is the Triangle Removal Lemma by Ruzsa and Szemer\'edi. To state it we use the asymptotic notation $o(.).$ For two functions over the reals, $f(x)$ and
$g(x),$ we write  $f(x)=o(g(x))$ if $f(x)/g(x)\rightarrow 0$ as $x\rightarrow \infty.$ The Triangle Removal Lemma states that any graph on $n$ vertices which contains at most $o(n^3)$  triangles can be make triangle free by removing at most $o(n^2)$ edges. See \cite{RUZSZE} for the original formulation of this result.

\begin{theorem}[Removal Lemma for Flats]\label{removal}
Given an arrangement of $n$  $d$-dimensional flats in $\R^{2d},$ such that no two flats are parallel. If the number of pairs with zero dimensional (single point) intersection is $o(n^2),$ then one can remove $o(n)$ flats such that the intersections of the remaining pairs have dimension at least one.
\end{theorem}

{\em Proof:} Let $G_n$ be a graph where the vertices represent the flats and two are connected iff the corresponding flats intersect in a single point. Let us suppose that the number of edges in $G_n$ is $\delta n^2$ for some $\delta>0,$ the minimum vertex cover\footnote{The vertex cover number of $G_n$ is the size of the smallest subset of vertices such that every edge has at least one endvertex in the set}  is  $m,$ and the maximum matching in $G_n$ has $M$ edges. Since for every maximal matching  the vertices form a vertex cover, we have $2M\geq m.$ By Lemma \ref{flat_pairs} we have the inequalities

\[ \delta n^2\geq\frac{M^2}{4^d}\geq\frac{m^2}{4^{d+1}},\]

implying

\[\sqrt{\delta}2^{d+1}n\geq m.\]

As  $\delta$ goes to zero, $m,$ the minimum number of flats needed to remove to avoid single vertex intersections, is getting smaller, it is $o(n)$ as we wanted to show.
\qed

\medskip

In the next corollary we show that an arrangement of two-dimensional flats has many zero-dimensional intersections in $\mathbb{R}^4$ unless there are some obvious obstacles, like many flats in a hyperplane, or flats intersecting in the same line. For two flats, $L$ and $F,$ the {\em affine span} or just {\em span} of them is the smallest dimensional flat which contains both $L$ and $F.$ 

\medskip

\begin{corollary}
Let $L_1, L_2, \ldots, L_n$ be 2-dimensional flats in $\mathbb{R}^4$ such that no three intersect in a single line and no two are parallel. If the number of $L_i,L_j$ $(1\leq i<j\leq n),$ pairs spanning $\mathbb{R}^4$ is $\delta n^2,$ then there is a hyperplane which contains  at least $(1-8\sqrt{\delta})n$ flats.
\end{corollary}

{\em Proof:} By Theorem \ref{removal} one can remove $8\sqrt{\delta}n$ flats so that the remaining flats intersect in a line. Let us select two flats from the remaining set, say $L_i$ and $L_j,$ they span a hyperplane in $\mathbb{R}^4.$ Any other flat, $E,$ is spanned by the two lines $L_i\cap E$ and $L_j\cap E,$ so $E$ is in the same hyperplane. (The two lines are distinct since no three flats intersect in a line.)
\qed

\section*{Acknowledgement}
The author is thankful to Frank de Zeeuw, Josh Zahl, and Istv\'an Tomon for the useful conversations.


\begin{thebibliography}{00}

\bibitem{AJF} S. Akbari, M. Jamaali, and S. A. Seyed Fakhari. {\em The clique numbers of regular graphs of matrix algebras are
finite}, Linear Algebra Appl., 2009., 431(10):1715--1718.

\bibitem{BT} A. Ben-Israel and T.N.E Greville,  {\em Generalized inverses: Theory and applications} (2nd ed.). 2003, New York, NY: Springer.


\bibitem{CLP} E. Croot, S. Lev, and P. Pach, {\em Progression-free sets in $\Z_4^n$ are exponentially small}, Ann. Math., 2017, 185(1):331--337.

\bibitem{EG} J. Ellenberg and D Gijswijt, {\em On large subsets of $\F_3^n$ with no three-term arithmetic
progression}, Ann. Math., 2017, 185 (1):339--343.

\bibitem{GUTH} L. Guth, {\em Polynomial Methods in Combinatorics}, 
AMS, University Lecture Series
Volume: 64; 2016; 273 pp;  

\bibitem{HENDERSON} H. V. Henderson and S. R. Searle,
{\em On Deriving the Inverse of a Sum of Matrices}
SIAM Review 1981 23:1, 53--60. 

\bibitem{MAMA} M. Marcus, {\em Determinants of Sums}, 
The College Mathematics Journal, 1990, Vol 21, No 2, 130--135.

\bibitem{RUZSZE} I.Z. Ruzsa and E. Szemer\'{e}di, {\em Triple systems with no six points carrying three triangles}, in Combinatorics (Keszthely, 1976), Coll. Math. Soc. J. Bolyai 18, Volume II, 939--945.



\bibitem{TAO} T. Tao, {\em Topics in Random Matrix Theory},
AMS, Graduate Studies in Mathematics
Volume: 132; 2012; 282 pp;  

\bibitem{TAO2} T. Tao, {\em The sum-product phenomenon in arbitrary rings}, Contrib. Discrete Math. 4 (2) (2009) 59--82.

\bibitem{TO} I. Tomon. {\em On the chromatic number of regular graphs of matrix algebras.} Linear Algebra Appl.,2015. 475:154--
162, 


\bibitem{TUR} P. Tur\'an,  {\em On an extremal problem in graph theory}, Matematikai \'es Fizikai Lapok (in Hungarian), 1941, 48: 436--452.

\end{thebibliography}
\end{document}